\documentclass[12pt]{article}
\usepackage[frenchb,english]{babel}
\usepackage[latin1]{inputenc}
\usepackage{amsthm}
\usepackage{amssymb,amsmath}
\newtheorem{prop}[subsubsection]{Proposition}

\newtheorem{surprop}[subsection]{Proposition} 
\newtheorem{surlem}[subsection]{Lemme} 
\newtheorem{surcor}[subsection]{Corollaire} 
\newtheorem{surthm}[subsection]{Th\'eor\`eme} 
\newtheorem{surconj}[subsection]{Conjecture} 
\newtheorem{lem}[subsubsection]{Lemme}

\theoremstyle{definition}
\newtheorem{para}[subsection]{}

\newcommand{\dem}{{\bf D\'emonstration}} 

\newcommand{\what}{\widehat} 
\newcommand{\ot}{\otimes}

\newcommand{\rig}{\to}

\newcommand{\hrig}{\hookrightarrow} 
\newcommand{\sta}{\stackrel}

\newcommand{\lrig}{\gets}

\newcommand{\der}{\partial} 
\newcommand{\la}{\langle} 
\newcommand{\ra}{\rangle}

\newcommand{\Qr}{{\bf Q}}
 
\newcommand{\Ne}{{\bf N}}

\newcommand{\Rr}{{\bf R}}

\newcommand{\Ga}{\Gamma}

\newcommand{Å}{{\cal A}} 
\newcommand{\BB}{{\cal B}} 
 
\newcommand{\DD}{{\cal D}} 
\newcommand{\EE}{{\cal E}} 
\newcommand{\FF}{{\cal F}}

\newcommand{\HH}{{\cal H}}

\newcommand{\OO}{{\cal O}}

\renewcommand{\SS}{{\cal S}} 
\newcommand{\TT}{{\cal T}} 
\newcommand{\UU}{{\cal U}} 
\newcommand{\VV}{{\cal V}} 
\newcommand{\XX}{{\cal X}} 
 
\newcommand{\ZZ}{{\cal Z}} 

\newcommand{\Ddag}{\DD^{\dagger}} 
\newcommand{\spec}{\rm spec} 
\begin{document} 
\date{}
\title{Sur l'holonomie de $\cal D$-modules arithm\'etiques associ\'es 
\`a des $F$-isocristaux surconvergents sur des courbes lisses} 
\author{Christine Noot-Huyghe et Fabien Trihan
\footnote{Pour cette collaboration, les deux auteurs ont
  b\'en\'efici\'e du soutien du Fond National de Recherche Scientifique (Belgique), ainsi que de l'IRMA, laboratoire 
de math\'ematiques de l'Universit\'e Louis Pasteur (Strasbourg). Ce travail a bénéficié en
outre du
soutien du réseau européen de recherche Arithmetic Algebraic Geometry de la Communauté
Européenne (Programme FP6, contrat MRTN-CT2003-504917).}} 
\maketitle 
\selectlanguage{frenchb}
\begin{abstract} Nous montrons que le $\cal D$-module arithmétique associé à un 
$F$-isocristal surconvergent sur une courbe lisse est holonome. Nous montrons 
d'abord que les $F$-isocristaux unipotents sont des $\cal D$-modules holonomes
en utilisant le fait que de tels $F$-isocristaux proviennent de $F$-isocristaux 
logarithmiques. Nous déduisons le cas général du théorème de réduction 
semi-stable pour les $F$-isocristaux sur les courbes de Matsuda-Trihan qui repose 
sur le théorème de monodromie $p$-adique démontré indépendamment par André, Kedlaya 
et Mebkhout.
\end{abstract}
\selectlanguage{english}
\begin{abstract}
\noindent We show that the arithmetic $\cal D$-module associated to an overconvergent $F$-isocrystal over 
a smooth curve is holonomic. We first prove that unipotent $F$-isocrystals are holonomic
${\cal D}$-module by using the fact that such $F$-isocrystals come from logarithmic
$F$-isocrystals. We deduce the general case from the semi-stable reduction theorem for $F$-isocrystals 
over curves of Matsuda-Trihan which relies on the $p$-adic monodromy theorem independently proved
by André, Kedlaya and Mebkhout.
\end{abstract}
\selectlanguage{frenchb}
\section*{Introduction} 
 Dans \cite{caro}, Caro montre que le 
 ${\cal D}$-module arithm\'etique associ\'e \`a un $F$-isocristal  
surconvergent sur une courbe lisse est holonome. Nous proposons de fournir une seconde preuve de ce r\'esultat. 
Notre d\'emonstration diff\`ere de celle de \cite{caro} sur les points suivants: La preuve de Caro repose sur 
une caract\'erisation 
des $F$-isocristaux surconvergents due \`a Berthelot. Nous proposons une preuve qui ne fait pas appel \`a ce r\'esultat. 

Notre preuve 
repose en effet sur des r\'esultats pr\'ec\'edents du second auteur concernant les $F$-isocristaux 
surconvergents unipotents. En particulier, il est d\'emontr\'e dans \cite{M-T_unip} que tout $F$-isocristal unipotent
 sur une courbe lisse provient d'un $F$-log cristal sur la compactification. Nous sommes alors en mesure de d\'ecrire 
le foncteur qui associe \`a tout $F$-log cristal un $F$-isocristal surconvergent sur le lieu o\`u la log structure
 devient triviale (cf \cite{LS-T_logcris}) de mani\`ere explicite en terme de ${\cal D}$-modules 
 arithm\'etiques. Le fait que le prolongement logarithmique de l'isocristal soit muni d'un Frobenius implique automatiquement
 l'\'egalit\'e entre la cohomologie de ces deux coefficients et permet de simplifier la preuve de \cite{caro}
 o\`u la structure de Frobenius sur le log cristal n'est pas utilis\'ee. 

Pour généraliser ce résultat aux $F$-isocristaux surconvergents, nous utilisons le
théorème de réduction semi-stable des $F$-isocristaux surconvergents dans le cas des
courbes dû à Matsuda-Trihan 
(\cite{M-T_unip}). Notons que ce résultat repose essentiellement sur le théorème de 
monodromie $p$-adique démontré indépendamment par Kedlaya (\cite{kedlaya_monp}), Mebkhout
(\cite{mebkhout_monp}), André (\cite{andre_monp}). Techniquement, nous généralisons au cas
des $\DD$-modules arithmétiques cohérents à pôles surconvergents, un résultat antérieur de 
Tsuzuki pour les images directes et inverses d'isocristaux surconvergents par un morphisme 
génériquement étale. En particulier, on construit un morphisme trace dans ce contexte.

Nous terminons par
 examiner la caract\'erisation des
 isocristaux surconvergents de Berthelot sous sa forme locale et globale et montrons 
comment
 cette derni\`ere permet de d\'emontrer la conjecture D de Berthelot,
 achevant ainsi la d\'emonstration de la stabilit\'e des $\cal 
D$-modules arithm\'etiques sur des courbes lisses par les six op\'erations
 cohomologiques de Grothendieck \`a l'aide des travaux de Caro
 (\cite{Car_FL}, \cite{caro}, \cite{caro2}, \cite{caro3}).\\

Nous tenons enfin \`a remercier Pierre Berthelot pour ses précisions et corrections concernant la partie
4, ainsi que le r\'ef\'er\'e pour ses judicieuses remarques.
\section{Notations-Généralités}
\subsection{Notations et Rappels}\label{subsection-notations}
Soit $k$ un corps parfait de caract\'eristique $p$, $W$ son  
anneau des vecteurs de Witt, $K= Frac(W)$, $\SS=Spf W$, 
$S_i=\spec W/p^{i+1}W$. Soient ${\cal X}/\SS$ un  
sch\'ema formel propre et lisse de dimension relative $N$, $\cal Z$ un  
diviseur \`a croisements normaux relatif de $\cal X$ et $\cal U:=X\backslash Z$.  
On notera respectivement $X$, $Z$ et $U$ leur r\'eduction modulo $p$. Ces 
sch\'emas sont naturellement munis d'une  log-structure: nous notons 
ainsi ${\cal X}^\#$ (resp. ${\cal Z}^\#$) le  log-sch\'ema de sch\'ema 
sous-ja\c cent $\cal X$ (resp. $\cal Z$) et dont la log-structure est 
induite par $\cal Z$ (resp. par la  log-structure inverse de celle 
de ${\cal X}^\#$) de sorte qu'on obtient un diagramme commutatif de 
log-sch\'emas 
$$\begin{array}{rcl}{\cal Z}^\# & \sta{\pi}{\to} & {\cal Z}\\ 
u^\#\downarrow& &\downarrow u \\
 {\cal X}^\# & \sta{\pi}{\to}&{\cal X}.\end{array}$$ 
Notons $X_i^\#/S_i$ la r\'eduction mod $p^{i+1}$ de ${\cal X}^\#/W$.  
 
L'immersion diagonale de $X_i^\#$ se factorise en 
$$X_i^\#\hookrightarrow X_i^\#(1)\to X_i^\#×_{S_i} X_i^\#$$ 
o\`u la premi\`ere fl\`eche est une immersion ferm\'ee exacte (cf \cite{Kato_log1}).  
Suivant ~\cite{Monta_th}, on note alors ${\cal P}^n_ 
{X_i^\#(m)}$ (resp. ${\cal P}^n_{X_i^\#}$) le faisceau des 
parties principales de niveau $m$ et d'ordre n de l'immersion 
fermée exacte $$X_i^\#\hookrightarrow X_i^\#(1)$$ 
(resp.  le faisceau des parties principales de niveau $m$ et d'ordre n de l'immersion 
fermée $X_i\hookrightarrow X_i(1)$). Ils sont tous deux munis d'une structure 
de ${\cal O}_{X_i}$-module. On note ${\cal D}_{X_i^\#,n}^{(m)}$ 
 (resp. ${\cal D}_{X_i^\#,n}$) son ${\cal O}_{X_i}$-module dual. Ce  
faisceau est appel\'e le faisceau des op\'erateurs diff\'erentiels de  
niveau $m$ et d'ordre $n$ (resp. le faisceau des op\'erateurs  
diff\'erentiels d'ordre n) sur $X_i^\#$. On note 
$${\cal  D}_{X_i^\#}^{(m)}:=\bigcup_n {\cal D}_{X_i^\#,n}^{(m)}\quad  
({\rm resp.}\, {\cal  D}_{X_i^\#}:=\bigcup_n {\cal D}_{X_i^\#,n}),$$  
$$\what{\cal D}^{(m)}_{{\cal X}^\#}:=\varprojlim {\cal  D}_{X_i^\#}^{(m)}\quad  
({\rm resp.}\, \what{\cal D}_{{\cal  X}^\#}:=\varprojlim {\cal D}_{X_i^\#})$$
 et finalement ${\cal  D}_{{\cal X}^\#}^\dagger:=\bigcup_m \what{\cal D}^{(m)}_{{\cal X}^\#}$  
($\subset \what{\cal D}_{{\cal X}^\#}$). On montre de mani\`ere 
analogue \`a \cite{Be1}, 3.6 la coh\'erence du faisceau ${\cal  D}_{{\cal 
X}^\#}^\dagger$. 

On dispose \'egalement des faisceaux  d'op\'erateurs 
diff\'erentiels classiques ${\cal D}_{{\cal  X}}^\dagger$ et 
$\what{\cal D}_{{\cal X}}$ d\'efinis dans \cite{Be1}. De plus,  
la fl\`eche $\pi:{\cal X}^\#\to {\cal X}$ induit un morphisme canonique injectif 
de faisceaux d'algèbres de $\OO_{\XX}$-algèbres 
$${\cal D}_{{\cal X}^\#}^\dagger\to {\cal D}_{{\cal X}}^\dagger.$$ 

Dans le cas d'un diviseur à croisements normaux, Montagnon montre dans le 
chapitre 5 de~\cite{Monta_th} que les faisceaux ${\cal  D}_{X_i^\#}^{(m)}$ sont de dimension 
cohomologique finie. En procédant comme en 4.4 de \cite{Be-smf}, on en déduit que les faisceaux 
$ \what{\cal D}^{(m)}_{{\cal X}^\#}$ et ${\cal  D}_{{\cal X}^\#,\Qr}^\dagger $ sont de dimension cohomologique finie.

Nous aurons aussi besoin du faisceau ${\cal D}_{{\cal X}}^\dagger(^{\dagger}Z)$, 
qui est une version à coefficients surconvergents 
du faisceau ${\cal D}_{{\cal X}}^\dagger$ et dont nous rappelons brièvement la 
définition. Soient $\VV$ un ouvert affine $\subset \XX$, $V_0$ la fibre spéciale de 
cet ouvert, $t$ une équation locale de $\ZZ$ sur $\VV$.  
On introduit le coefficient $$\BB^{(m)}_{X_i}=\OO_{X_i}(T)\left/t^{p^{m+1}}T-p\right., $$
puis $$\DD^{(m)}_{X_i}(Z)= \BB^{(m)}_{X_i}\ot_{\OO_{X_i}}\DD^{(m)}_{X_i},$$
et les faisceaux complétés $\what{\BB}^{(m)}_{\XX}$ et $\what{\DD}^{(m)}_{\XX}(Z)$.
Passons à la limite inductive sur $m$: on définit 
$$\OO_{\XX,\Qr}(^{\dagger}Z)=\varinjlim_m \what{\BB}^{(m)}_{\XX,\Qr} \quad{\rm et}\quad  
\Ddag_{\XX,\Qr}(^{\dagger}Z)=\varinjlim_m \what{\DD}^{(m)}_{\XX,\Qr}(Z) ,$$
o\`u si $\FF$ est une faisceau quelconque, on notera $\FF_\Qr$
pour $\FF\ot\Qr$. Les faisceaux considérés sont des faisceaux d'algèbres cohérents. On voit facilement 
que les faisceaux $\OO_{\XX,\Qr}(^{\dagger}Z)$ et $\Ddag_{\XX,\Qr}(^{\dagger}Z)$ ne dépendent 
que de la réduction mod $p$ du diviseur $\ZZ$ notée $Z$ et du sous-schéma réduit
$Z_{red}$, ce qui justifie la notation. 

Si $\FF$ est un faisceau de $\OO_{\XX,\Qr}$-modules, on pose 
$$\tilde{\FF}=\OO_{\XX}(^{\dagger }Z)\ot_{\OO_{\XX}}\FF.$$

On dispose de morphismes canoniques
$$ \OO_{\XX,\Qr}\hrig \OO_{\XX,\Qr}(^{\dagger}Z) \quad{\rm et}\quad
{\cal D}_{{\cal X},\Qr}^\dagger\hrig {\cal D}_{{\cal X},\Qr}^\dagger(^{\dagger}Z).$$
  
Nous noterons \'egalement ${\bf res}$, le foncteur qui \`a tout 
${\cal D}^\dagger_{{\cal X},\Qr}(^\dagger Z)$-module associe sa structure canonique de 
${\cal D}^\dagger_{{\cal X},\Qr}$-module via le morphisme canonique précédent.

Nous adopterons le formalisme des images directes et inverses au sens des ${\cal D}$-modules 
arithmétiques introduits dans \cite{Be-survey}. Soient $\XX$ 
et $\XX'$ des schémas formels lisses sur $\SS$, $f$: $\XX'\rig \XX$ un morphisme, 
$d=dim\XX'-dim\XX$, 
$\UU'\subset \XX'$ le complémentaire d'un diviseur $\ZZ'$, $\UU$ le complémentaire d'un diviseur 
$\ZZ$ de $\XX$, tel que $f(\UU')\subset \UU$. On peut supposer pour les 
constructions que l'on a 
un morphisme de faisceaux d'algèbres 
$f^{-1}\BB^{(m)}_{\XX}\rig \BB^{(m)}_{\XX'}$, qui induit 
(1.1.3 de~\cite{caro}) pour tout $i$ des morphismes
$f^{-1}\BB^{(m)}_{X_i}\rig \BB^{(m)}_{X'_i}$. 
On construit à partir de là les faisceaux 
$$\DD^{(m)}_{X_i'\rig X_i}(Z')=\BB^{(m)}_{X_i'}\ot_{\OO_{X_i'}}\DD^{(m)}_{X_i'\rig X_i}, 
\what{\DD}^{(m)}_{\XX'\rig \XX}(Z')=\varprojlim_i\DD^{(m)}_{X_i'\rig X_i}(Z'),$$
où $\DD^{(m)}_{X_i'\rig X_i}$ est introduit en 2.2 de \cite{Be-survey}, et 
$$\DD^{(m)}_{X_i\lrig X_i'}(Z')=\BB^{(m)}_{X_i'}\ot_{\OO_{X_i'}}\DD^{(m)}_{X_i\lrig X_i'}, 
\quad \what{\DD}^{(m)}_{\XX\lrig \XX'}(Z')=\varprojlim_i\DD^{(m)}_{X_i\lrig X_i'}(Z'),$$
où $\DD^{(m)}_{X_i\lrig X_i'}$ est introduit en 2.2 de \cite{Be-survey}. 
On pose aussi 
$$\Ddag_{\XX'\rig \XX,\Qr}(^{\dagger}Z')= \varinjlim_m\what{\DD}^{(m)}_{\XX'\rig \XX,\Qr}(Z'), \quad 
\Ddag_{\XX\lrig \XX',\Qr}(^{\dagger}Z')= \varinjlim_m\what{\DD}^{(m)}_{\XX\lrig \XX',\Qr}(Z').$$
Soit $\EE$ un faisceau de $\Ddag_{\XX,\Qr}(^{\dagger}Z)$-modules cohérents, on définit 
 $\tilde{f}^!\EE$ comme dans 4.3.3 de \cite{Be-survey}. Soit 
$\EE'$ un faisceau de $\Ddag_{\XX',\Qr}(^{\dagger}Z')$-modules cohérents, on définit 
 $\tilde{f}_+\EE'$ comme dans 4.3.6 de \cite{Be-survey}. Les versions sans diviseur
(i.e $Z$ et $Z'$ sont vides) de ces foncteurs sont respectivement notés 
$f^!\EE'$ et $f_+\EE$. 

\vspace{+8mm}
Pour fixer les idées, donnons des descriptions en coordonnée locale, sur une courbe, des faisceaux
intervenants, sous nos hypothèses que $Z$ est un diviseur à croisements normaux de $X$.
\subsection{Descriptions en coordonnée locale}
\label{subsection-desc_loc}
Soit $\XX$ une courbe. Plaçons-nous sur un ouvert $\VV$ de $\XX$, muni d'une coordonnée $t$, telle que 
$\ZZ\bigcap\VV=V(t)$. Les symboles $a_k$ et $a_{l,k}$ qui suivent 
désignent des éléments de $\Ga(\VV,\OO_{\XX,\Qr})$, que l'on munit de la norme 
spectrale. On note $\partial $ la dérivée partielle par rapport à $t$, et, comme 
d'habitude, $\partial^{[k]}=\partial^k/k!$. Nous avons alors les descriptions suivantes des faisceaux 
introduits ci-dessus
$$\Ga(\VV,\Ddag_{\XX,\Qr})=\left\{\sum_{k\in\Ne}a_k\partial^{[k]}\,|\,\exists C>0, \eta<1 \, {\rm tels} \,
{\rm que} \, |a_k|_{sp}<C\eta^k \right\},$$
$$\Ga(\VV,\Ddag_{\XX,\Qr}(^{\dagger} Z))=\left\{\sum_{l,k\in\Ne}a_{l,k}t^{-l}\partial^{[k]}\,|\,\exists C>0, \eta<1 \, {\rm tels} \,
{\rm que} \, |a_{l,k}|_{sp}<C\eta^{l+k} \right\}.$$
Notons enfin  
$$\der_{[ k ]}=\frac{1}{k!}(t\der)(t\der-1)(t\der-2)\cdots
(t\der-k+1).$$ Alors, d'après 2.3 C de \cite{Monta_th}, on a 
$$\Ga(\VV,\Ddag_{\XX^{\#},\Qr})=\left\{\sum_{k\in\Ne}a_k\der_{[ k ]}\,|\,\exists C>0, \eta<1 \, {\rm tels} \,
{\rm que} \, |a_k|_{sp}<C\eta^k \right\}.$$
\subsection{Morphisme trace pour les $\DD$-modules arithmétiques à coefficients
surconvergents}
On considère ici $\XX$ 
et $\XX'$ des schémas formels lisses sur $\SS$, $f$: $\XX'\rig \XX$ un morphisme fini, 
de degré générique égal à $d$, 
$\UU'=f^{-1}(\UU)$, $\UU$ le complémentaire d'un diviseur 
$\ZZ$ de $\XX$, tel qu'on ait un diagramme commutatif 
$$ \begin{array}{rcl}{\cal U}'&\sta{j'}{\rig}&{\cal X'}\\
 g\downarrow&&\downarrow f \\ 
{\cal U}&\sta{j}{\rig}&{\cal X},\end{array}$$ 
où la restriction $g$ de $f$ à $\UU'$ est finie et étale. On considère 
la réduction mod $p$ de ces schémas et de ce diagramme ainsi que 
$f_0$ l'homomorphisme induit par $f$: $X'\rig X$. Dans un contexte analogue, 
Tsuzuki construit l'image directe et inverse par $f_0$ d'isocristaux surconvergents sur 
$U$ et $U'$. On donne ici des résultats analogues 
pour l'image inverse par $f$ des  
$\Ddag_{\XX,\Qr}(^{\dagger}Z)$-modules cohérents et l'image directe des $\Ddag_{\XX',\Qr}(^{\dagger}Z')$-modules 
cohérents. Comme $\ZZ'=f^{-1}(\ZZ)$, on peut supposer que
$\what{\BB}^{(m)}_{\XX'}=f^*\what{\BB}^{(m)}_{\XX}$, pour faire les constructions d'images
directes (resp. inverses) des $\Ddag_{\XX',\Qr}(^{\dagger}Z')$-modules (resp.
des $\Ddag_{\XX,\Qr}(^{\dagger}Z)$ modules) comme en 7 de~\cite{HuFour}. Dans la section 
7 de~\cite{HuFour}, on considère la cas particulier où $g$ est un isomorphisme. On
explique ici ce qu'il faut ajouter aux considérations de~\cite{HuFour} pour obtenir le 
cas considéré ici. 
Partons du lemme crucial suivant (7.2.3 de~\cite{HuFour}).
\begin{lem}\label{subsubsection-lemcruc} \begin{enumerate}
\item[(i)] Il existe un isomorphisme canonique de 
$\Ddag_{\XX',\Qr}(^{\dagger}Z')$-modules à gauche
$$\Ddag_{\XX',\Qr}(^{\dagger}Z')\sta{\sim}{\rig}\Ddag_{\XX' \rig \XX,\Qr}(^{\dagger}Z').$$
On obtient ainsi un morphisme canonique injectif de faisceaux de $f^{-1}\OO_{\XX,\Qr} $-algèbres 
$$f^{-1}\Ddag_{\XX,\Qr}(^{\dagger}Z)\hrig \Ddag_{\XX',\Qr}(^{\dagger}Z').$$
\item[(ii)] Il existe un isomorphisme canonique de 
$\Ddag_{\XX',\Qr}(^{\dagger}Z')$-modules à droite
$$\Ddag_{\XX',\Qr}(^{\dagger}Z')\sta{\sim}{\rig}\Ddag_{\XX \lrig\XX',\Qr}(^{\dagger}Z').$$
\end{enumerate}
\end{lem}  

\dem. La première partie de (i) est exactement 7.2.3 de~\cite{HuFour}, avec des hypothèses plus faibles. Nous
expliquerons ci-dessous pourquoi l'énoncé est vrai dans notre cas.
 
Justifions la deuxième partie du (i). 
On dispose d'un morphisme $f^{-1}\Ddag_{\XX,\Qr}(^{\dagger}Z)$-linéaire à droite
$f^{-1}\Ddag_{\XX,\Qr}(^{\dagger}Z)\rig f^*\Ddag_{\XX',\Qr}(^{\dagger}Z')$, 
qui envoie $P$ sur $1\ot P$. Grâce au (i) du lemme, cela donne un morphisme 
$\Ddag_{\XX,\Qr}(^{\dagger}Z)$-linéaire à droite 
de $\Ddag_{\XX,\Qr}(^{\dagger}Z)$ vers $\Ddag_{\XX',\Qr}(^{\dagger}Z')$. Il suffit de voir en
restriction à l'ouvert $\UU$ que c'est un morphisme de faisceaux d'algèbres, ce qui est 
clair car dans ce cas le morphisme considéré est un isomorphisme de faisceaux d'algèbres. 
Cela complète la première assertion qui ne figure pas stricto sensu dans~\cite{HuFour}. 

Pour le (ii), on procède comme précédemment et on utilise l'isomorphisme de transposition qui suit et qui 
identifie les deux faisceaux de $\OO_{\XX',\Qr}$-algèbres (1.3.4 de~\cite{Be-smf})
$$\Ddag_{\XX',\Qr}(^{\dagger}Z')\simeq \tilde{\omega}_{\XX',\Qr}\ot_{\OO_{\XX',\Qr}(^{\dagger} Z')}
\Ddag_{\XX',\Qr}(^{\dagger}Z')\ot_{\OO_{\XX',\Qr}(^{\dagger} Z')}
\tilde{\omega}_{\XX',\Qr}^{-1},$$
où $\omega_{\XX'}$ est le faisceau inversible des différentielles de rang maximal sur
$\XX'$.

Résumons la construction de la partie 7 de \cite{HuFour}. En réalité, on suppose 
dans cette partie 7 de \cite{HuFour} que $f$ induit un isomorphisme 
$f^{-1}\UU\simeq \UU$, mais on utilise seulement le fait qu'un certain 
déterminant jacobien est inversible en restriction à $\UU'$, de sorte  
que seule l'hypothèse que $f_{|\UU'}$ est étale intervient. Plus précisément, on constate, 
grâce à un lemme de géométrie rigide, que pour $m$ assez grand on a
un isomorphisme d'espaces tangents 
$$J\colon \what{\BB}^{(m)}_{\XX',\Qr}\ot_{\OO_{\XX'}} \TT_{\XX'}\sta{\sim}{\rig}
f^*\what{\BB}^{(m)}_{\XX,\Qr}\ot_{\OO_{\XX}}\TT_{\XX}.$$
Soit $j_m$ l'homomorphisme canonique 
$$j_m \colon \what{\DD}^{(m)}_{\XX',\Qr}(Z')\rig
f^*\what{\DD}^{(m)}_{\XX,\Qr}(Z).$$
Comme ce morphisme devient un isomorphisme en restriction à $\UU$, 
$j_m$ est injectif. 

On utilise $J$ pour construire, à partir d'un $m$ assez grand, une suite croissante 
d'entiers $u_m$ ($u_m\geq m$) et des injections continues pour la topologie 
$p$-adique 
$$i_m\colon f^*\what{\DD}^{(m)}_{\XX,\Qr}(Z){\hrig}
\what{\DD}^{(u_m)}_{\XX',\Qr}(Z'),$$
telles que $i_m\circ j_m$ est l'injection canonique 
$$\what{\DD}^{(m)}_{\XX',\Qr}( Z')\hrig \what{\DD}^{(u_m)}_{\XX',\Qr}( Z') .$$
 De plus, si $\der$ est une section locale de $\TT_{\XX}$,  
 alors $i_m(1\ot \der)=J^{-1}(1\ot\der)$.

Remarquons que, avec nos hypothèses ($\UU'=f^{-1}\UU$ et $f$ est fini), $\tilde{f}_+ $ 
préserve la cohérence (4.3.8 de ~\cite{Be-survey}).

\vspace{+4mm}
Le lemme permet de donner la description suivante de l'image directe d'un 
$\Ddag_{\XX',\Qr}(^{\dagger}Z')$-module cohérent (resp. de l'image inverse d'un 
$\Ddag_{\XX,\Qr}(^{\dagger}Z)$-module cohérent).
\begin{prop} ~\label{subsubsection-op_coh_gen_et}
\begin{enumerate}\item[(i)] On a un isomorphisme canonique de 
$\Ddag_{\XX}(^{\dagger}Z)$-modules à droite
$$\Ddag_{\XX'\rig\XX}(^{\dagger}Z')\simeq f^*\Ddag_{\XX}(^{\dagger}Z).$$
Soit $\EE$ un $\Ddag_{\XX,\Qr}(^{\dagger}Z)$-module cohérent. Alors il
existe un isomorphisme canonique de $\Ddag_{\XX',\Qr}(^{\dagger}Z')$-modules à gauche 
$$ \tilde{f}^!\EE\sta{\sim}{\rig}\Ddag_{\XX'}(^{\dagger}Z')
\ot_{f^{-1}\Ddag_{\XX}(^{\dagger}Z)}f^{-1}\EE.$$ 
Le module $\tilde{f}^!\EE$ est concentré en degré $0$ et
 est isomorphe comme $\OO_{\XX',\Qr}$-module à $f^*\EE=\OO_{\XX'}\ot_{f^{-1}\OO_{\XX}}f^{-1}\EE$. 
\item [(ii)] Soit $\EE'$ un $\Ddag_{\XX',\Qr}(^{\dagger}Z')$-module cohérent, alors 
$\tilde{f}_+\EE'\simeq f_*\EE'$ comme faisceau de
$\OO_{\XX,\Qr}(^{\dagger}Z)$-modules.\end{enumerate}
\end{prop} 
On peut bien entendu donner une variante de cette proposition pour 
$\EE\in D^b_{coh}(\Ddag_{\XX,\Qr}(^{\dagger}Z))$ et $\EE'\in
D^b_{coh}(\Ddag_{\XX',\Qr}(^{\dagger}Z'))$. 

\vspace{+2mm}
\dem. Par définition, on a  
$$f^{!}\EE=\Ddag_{\XX' \rig \XX}(^{\dagger}Z')\ot_{f^{-1}\Ddag_{\XX}(^{\dagger}Z)}^{\bf L}
f^{-1}\EE.$$ D'après le lemme précédent, on a donc
$$f^!\EE\sta{\sim}{\rig}\Ddag_{\XX'}(^{\dagger}Z')\ot_{f^{-1}\Ddag_{\XX}(^{\dagger}Z)}^{\bf L}
f^{-1}\EE.$$ On est donc ramené à montrer que $\Ddag_{\XX',\Qr}(^{\dagger}Z')$ est plat comme 
$f^{-1}\Ddag_{\XX,\Qr}(^{\dagger}Z)$-module à droite, et donc que 
$\what{\DD}^{(m)}_{\XX' \rig \XX,\Qr}(Z')$ est plat à droite sur 
$\what{\DD}^{(m)}_{\XX,\Qr}(Z)$ pour $m$ assez grand. Dans notre situation, on a 
$\what{\BB}^{(m)}_{\XX'}=f^*\what{\BB}^{(m)}_{\XX}$. La question de la platitude est locale sur 
$\XX$ qu'on suppose donc affine et muni de coordonnées locales. Dans ce cas, 
l'algèbre $\what{\BB}^{(m)}_{\XX',\Qr}(\XX')$ est finie sur $\what{\BB}^{(m)}_{\XX,\Qr}(\XX)$ 
et $\what{\DD}^{(m)}_{\XX,\Qr}(Z)$ est une algèbre de Banach pour la norme associée à 
la topologie $p$-adique. Ainsi, le module 
$$\what{\BB}^{(m)}_{\XX',\Qr}(\XX')\ot_{\what{\BB}^{(m)}_{\XX,\Qr}(\XX)}
\what{\DD}^{(m)}_{\XX}(Z),$$
est de type fini sur $\what{\DD}^{(m)}_{\XX,\Qr}(Z)$. Ce module est donc 
complet pour la norme quotient provenant de sa structure de module et qui définit 
la topologie $p$-adique sur ce module. Ainsi, ce module s'identifie à 
$\what{\DD}^{(m)}_{\XX' \rig \XX,\Qr}(Z')$. Pour $m$ assez grand, 
le morphisme $f^{-1}\what{\BB}^{(m)}_{\XX,\Qr}\rig \what{\BB}^{(m)}_{\XX',\Qr}$ 
est fini et étale donc plat. 
On en déduit, pour $\XX$ général, un isomorphisme canonique de
$f^{-1}\what{\DD}^{(m)}_{\XX,\Qr}(Z)$-modules à droite 
$$\what{\DD}^{(m)}_{\XX'\rig\XX,\Qr}(Z')\simeq \what{\BB}^{(m)}_{\XX',\Qr}
\ot_{f^{-1}\what{\BB}^{(m)}_{\XX,\Qr}}f^{-1}\what{\DD}^{(m)}_{\XX}(Z).$$
Cela donne que $\what{\DD}^{(m)}_{\XX'\rig\XX,\Qr}(Z')$ est plat sur 
$f^{-1}\what{\DD}^{(m)}_{\XX,\Qr}(Z)$ et donc la platitude cherchée. 
De plus, comme $\what{\BB}^{(m)}_{\XX',\Qr}\simeq f^*\what{\BB}^{(m)}_{\XX,\Qr}$,  
 on en déduit l'isomorphisme 
$$\what{\DD}^{(m)}_{\XX'\rig\XX,\Qr}(Z')\simeq \OO_{\XX',\Qr}\ot_{f^{-1}\OO_{\XX,\Qr}}
f^{-1}\what{\DD}^{(m)}_{\XX}(Z).$$
En passant à la limite inductive sur $m$, on trouve l'isomorphisme 
$\Ddag_{\XX,\Qr}(^{\dagger}Z)$-linéaire à droite
$$\Ddag_{\XX'\rig \XX,\Qr}(^{\dagger}Z')\simeq f^*\Ddag_{\XX,\Qr}(^{\dagger}Z),$$
d'où l'ensemble des énoncés du (i).

En passant à la limite inductive 
sur $m$, et la deuxième partie du (i).
Au passage, on vérifie que le faisceau d'algèbres $\OO_{\XX'}(^{\dagger}Z')$ est fini 
et étale (donc plat) sur $\OO_{\XX}(^{\dagger}Z)$.

Comme $f$ est un morphisme fini, $f_*={\bf R}f_*$. 
Par définition, on a
$$f_+\EE'={\bf R}f_*\left(\Ddag_{\XX \lrig \XX'}(^{\dagger}Z')
\ot^{{\bf L}}_{\Ddag_{\XX'}(^{\dagger}Z')}\EE'\right).$$
L'énoncé (ii) suit donc directement du lemme précédent.

Précisons l'action des opérateurs différentiels par ces opérations cohomologiques. 
Soient $\EE$ un $\Ddag_{\XX,\Qr}(^{\dagger}Z)$-module cohérent, $e $ une
section locale de $\EE$, 
 $\der'$ une section locale du faisceau 
tangent $\TT_{\XX'}$. 
Alors, on a la description suivante de l'action de $\Ddag_{\XX',\Qr}(^{\dagger}Z')$ 
sur $f^*\EE$
$$\der'\cdot(1\ot e)=J(1\ot\der')\cdot (1\ot e) .$$
Soient $\EE'$ un $\Ddag_{\XX',\Qr}(^{\dagger}Z')$-module cohérent, 
$\VV$ un ouvert affine de $\XX$, $\der$ une section 
 du faisceau tangent $\TT_{\XX}$ sur $\VV$, $e'$ une section de
$f_*\EE'(\VV)$. 
Alors, on a la description suivante de l'action de $\Ddag_{\XX,\Qr}(^{\dagger}Z)$ 
sur $f_*\EE$
$$\der\cdot e'=J^{-1}(1\ot\der)\cdot e'\,\in f_*(\EE)(\VV) .$$

Ces résultats nous permettent de construire un morphisme trace dans ce contexte.
Commençons par interpréter l'homomorphisme d'adjonction dans la catégorie des 
$\Ddag_{\XX,\Qr}(^{\dagger}Z)$-modules cohérents.
\begin{prop} Soit $\EE$ un $\Ddag_{\XX,\Qr}(^{\dagger}Z)$-module cohérent. L'homomorphisme 
d'adjonction canonique $ad$ : $\EE\rig f_*f^*\EE$ est $\Ddag_{\XX,\Qr}(^{\dagger}Z)$-linéaire si 
on identifie $f_*f^*\EE$ à ${\tilde{f}}_+{\tilde{f}}^!\EE$. 
\end{prop} 
\dem. On a par définition $ad(e)=1\ot e$. L'assertion est locale sur la base. On peut 
donc supposer que $\XX$ est affine, muni de coordonnées locales et que 
$\ZZ=V(h)$ avec $h\in \OO_{\XX}(\XX)$. Comme ${\tilde{f}}^! $ est ${\tilde{f}}_+$ sont exacts, toute 
surjection $\Ddag_{\XX,\Qr}(^{\dagger}Z)$-linéaire $\Ddag_{\XX,\Qr}(^{\dagger}Z)^a\rig \EE$ 
donne lieu à une surjection $\Ddag_{\XX}(^{\dagger}Z)$-linéaire 
$f_*f^*\Ddag_{\XX}(^{\dagger}Z)\rig f_*f^*\EE$. 
Il suffit finalement de vérifier que 
l'adjonction est $\Ddag_{\XX,\Qr}(^{\dagger}Z)$-linéaire pour 
$\EE=\Ddag_{\XX,\Qr}(^{\dagger}Z)$. Or, dans ce cas, l'adjonction correspond 
à l'homomorphisme de faisceaux d'algèbres 
$\Ddag_{\XX,\Qr}(^{\dagger}Z)\hrig f_*(\Ddag_{\XX',\Qr}(^{\dagger}Z'))$. Cela donne la linéarité de 
l'adjonction. 
\begin{prop} \label{subsubsection-trace}
Soit $\EE$ un $\Ddag_{\XX,\Qr}(^{\dagger}Z)$-module cohérent. Alors il existe 
un morphisme trace $tr$: $ f_*f^*\EE\rig \EE$ qui est $\Ddag_{\XX,\Qr}(^{\dagger}Z)$-linéaire 
si on identifie $f_*f^*\EE$ à ${\tilde{f}}_+{\tilde{f}}^!\EE$. De plus, l'homomorphisme composé
$$\EE\sta{ad}{\rig}{\tilde{f}}_+{\tilde{f}}^!\EE \sta{tr}{\rig}\EE,$$
est égal à $d\cdot Id_{\EE}$. 
\end{prop} 
\dem. Pour $m$ assez grand, on a un morphisme fini étale 
$f^{-1}\what{\BB}^{(m)}_{\XX,\Qr}\rig \what{\BB}^{(m)}_{\XX',\Qr}$ qui permet 
de définir un morphisme trace $\what{\BB}^{(m)}_{\XX,\Qr}$-linéaire
$$t_m\colon f_*\what{\BB}^{(m)}_{\XX',\Qr}\rig \what{\BB}^{(m)}_{\XX,\Qr}.$$
Identifions 
$$f_*\OO_{\XX',\Qr}(^{\dagger}Z')\simeq
f_*\what{\BB}^{(m)}_{\XX',\Qr}\ot_{\what{\BB}^{(m)}_{\XX,\Qr}}\OO_{\XX}(^{\dagger}Z).$$
En étendant les scalaires, on dispose ainsi d'un morphisme trace 
$\OO_{\XX,\Qr}(^{\dagger}Z)$-linéaire $t=t_m\ot Id_{\OO_{\XX}(^{\dagger}Z)}$: 
$f_*\OO_{\XX',\Qr}(^{\dagger}Z')\rig \OO_{\XX,\Qr}(^{\dagger}Z).$ Soit $\EE$ un 
$\Ddag_{\XX,\Qr}(^{\dagger}Z)$-module cohérent. Comme $f$ est fini, 
on peut identifier $$f_*\left(\OO_{\XX'}(^{\dagger}Z')\ot_{f^{-1}\OO_{\XX}
(^{\dagger}Z)}f^{-1}\EE\right)\simeq f_*\left(\OO_{\XX'}(^{\dagger}Z')\right)
\ot_{\OO_{\XX}(^{\dagger}Z)}\EE.$$ 
Le morphisme $t\ot Id_{\EE}$ 
définit alors un morphisme $\OO_{\XX,\Qr}(^{\dagger}Z)$-linéaire 
$ f_*f^*\EE \rig  \EE .$
Il reste à vérifier que la flèche obtenue est $\Ddag_{\XX,\Qr}(^{\dagger}Z)$-linéaire. Par le 
même argument que précédemment, on se ramène au cas où $\EE=\Ddag_{\XX,\Qr}(^{\dagger}Z)$. 
Dans ce cas, la structure de $\Ddag_{\XX,\Qr}(^{\dagger}Z)$-module de 
$f_*(\Ddag_{\XX',\Qr}(^{\dagger}Z'))$ est donnée par l'homomorphisme de faisceaux d'algèbres
$\Ddag_{\XX}(^{\dagger}Z)\hrig f_*(\Ddag_{\XX'}(^{\dagger}Z'))$ et la trace qui est égale 
à $t\ot Id_{\Ddag_{\XX,\Qr}(^{\dagger}Z)}$ est $\Ddag_{\XX,\Qr}(^{\dagger}Z)$-linéaire. 
\subsection{Quelques résultats de compatibilité.}
On considère dans cette partie un schéma formel $\XX$ lisse sur $\SS$, muni d'un diviseur 
relatif $\ZZ$, ainsi que les faisceaux $\Ddag_{\XX}$ et $\Ddag_{\XX}(^{\dagger}Z)$.
On montre ici des résultats de compatibilité entre l'image inverse et l'image directe au 
sens des $\Ddag$-modules et au sens des $F$-isocristaux surconvergents. On commence par 
l'image inverse et un énoncé (1.5.4 de \cite{huy}) que nous reproduisons ici car il n'a pas été 
publié et est souvent cité. Soient $\XX$ 
et $\XX'$ des schémas formels lisses sur $\SS$, $f$: $\XX'\rig \XX$ un morphisme, 
$d=dim\XX'-dim\XX$, 
$\UU'\subset \XX'$ le complémentaire d'un diviseur $\ZZ'$, $\UU$ le complémentaire d'un diviseur 
$\ZZ$ de $\XX$, tel que $f(\UU')\subset \UU$. On note $f_0$ le morphisme induit
par $f$ au niveau des fibres spéciales, $f_0$ : $X'\rig X$.

Soit $E$ un isocristal surconvergent sur $U$ le long de $X\backslash Z$.  
Son image directe par spécialisation $sp_*E$ sur $\XX$ est un 
$\Ddag_{\XX}(^{\dagger} Z)$-module cohérent. On note $f_0^{rig,*}E$ l'image inverse de 
$E$ comme isocristal surconvergent sur $U'$.
On a alors le théorème de comparaison suivant
\begin{prop}\label{subsubsection-comp_iminv}
Il existe un isomorphisme canonique de $\Ddag_{\XX',\Qr}(^{\dagger} Z')$-modules 
à gauche $$sp_*f_0^{rig,*}E\sta{\sim}{\rig}\tilde{f}^!(sp_*E)[-d].$$
\end{prop} 
En particulier, $\tilde{f}^!(sp_*E)[-d]$ est un $\Ddag_{\XX'}(^{\dagger} Z')$-module cohérent.\\

Terminons maintenant par un énoncé de comparaison des images directes en tant 
que $\Ddag$-module et en tant qu'isocristal surconvergent sous les hypothèses de
 la sous-section précédente, i.e. $f$ est fini, génériquement étale sur 
un ouvert $\UU$. On consultera 5. de \cite{Ts1} pour les descriptions de l'image directe et inverse par 
$f_0$ des isocristaux surconvergents dans ce contexte.
\begin{prop}Soit $E'$ un $F$-isocristal surconvergent sur $U'$ (le long de $Z'$). Alors on 
a un isomorphisme canonique $$sp_*f_{0,*}^{rig}E'\sta{\sim}{\rig}\tilde{f}_+(sp_*E').$$
\end{prop} 
\dem. Posons $\EE'=sp_*E'$. C'est un $\Ddag_{\XX',\Qr}(^{\dagger}Z')$-module cohérent. 
Il résulte de \cite{Ts1}, 5.1.2 et de~\ref{subsubsection-op_coh_gen_et} que 
ces deux modules s'identifient à $f_*\EE'$ comme $\OO_{\XX,\Qr}(^{\dagger}Z)$-modules. 
Nous sommes donc ramenés à montrer que l'action de $\Ddag_{\XX,\Qr}(^{\dagger}Z)$ est
identique sur ces deux modules. La question est locale sur $\XX$, qu'on peut supposer 
affine, lisse, muni de coordonnées locales $x_1,\ldots,x_N$, et tel que 
$\ZZ=V(h)$ avec $h\in \OO_{\XX}(\XX)$. Comme $f_*\EE'$ est un
$\OO_{\XX,\Qr}(^{\dagger}Z)$-module localement projectif, on a une 
injection $f_*\EE'(\XX)\subset f_*\EE'(\UU)$. On a une injection analogue pour 
$\Ddag_{\XX,\Qr}(^{\dagger}Z)$. Il suffit donc de vérifier que l'action de 
$\Ddag_{\XX,\Qr}(^{\dagger}Z)$ est la même sur $sp_*f_{0,*}^{rig}E'$ et 
sur $\tilde{f}_+\EE'$, en restriction à $\UU$. Comme $f$ est fini et étale au-dessus de 
$\UU$, l'action de $\DD^{(0)}_{\UU,\Qr}$ est obtenue sur chacun de ces deux modules en prenant 
l'image directe de la $PD$-stratification définissant la connexion sur $\EE'$. Ainsi 
l'action des dérivations est identique sur $sp_*f_{0,*}^{rig}E'(\UU)$ et sur
$\tilde{f}_+\EE'(\UU)$, qui s'identifient tous les deux à 
$f_*\EE'(\UU)$. L'action des opérateurs 
$$\der_{x_i}^{{\la k_i \ra}_{(m)}}=\frac{q_{k_i}!}{k_i!}\der_{x_i}^{k_i},$$ 
où $q_{k_i}$ est le quotient de la division euclidienne de $k_i$ par 
$p^m$, est donc 
aussi identique. Comme $f_*(\EE')(\UU)$ est $p$-adiquement complet, cela 
donne que l'action de $\what{\DD}^{(m)}_{\UU,\Qr}(\UU)$ est la même et en passant 
à la limite inductive sur $m$, on voit que l'action de 
$\Ddag_{\UU,\Qr}$ est la même sur $sp_*f_{0,*}^{rig}E'$ et sur 
$\tilde{f}_+\EE'$ au-dessus de $\UU$. Cela montre la proposition.

\vspace{+8mm}
Dans toute la suite de cet article, nous reprenons les notations de
~\ref{subsection-notations} et supposons de plus que $\XX$ est de dimension relative 
$1$ sur $\SS$. 

\section{Holonomie des F-isocristaux unipotents} 
Nous montrons tout d'abord comment d\'emontrer l'holonomie des modules  
diff\'erentiels associ\'es \`a des $F$-isocristaux surconvergents unipotents. Si 
${\cal A}$ est un faisceau de $\OO_{\XX}$-algèbres, on note comme d'habitude 
$D^*_{coh}({\cal A})$ pour $*\in\{b,-\}$ la catégorie dérivée des complexes de 
${\cal A}$-modules  
à cohomologie bornée si $*=b$ et concentrée en degrés négatifs si $*=-$.
\subsection{Images directes de $\DD$-modules logarithmiques} 
Nous pouvons aussi construire de mani\`ere analogue \`a \cite{Be-survey},4.3.7.1 
l'image directe 
$$\pi_+:D^b_{coh}({\cal D}^\dagger_{{\cal X}^\#,\Qr})\to D^-({\cal D}^\dagger_{{\cal X},\Qr})$$ 
 $${\cal E}^.\mapsto {\cal D}^\dagger_{{\cal X}\gets {\cal X}^\#,
{\bf Q}}\otimes^{\bf L}_{{\cal D}^\dagger_{{\cal X}^\#}} {\cal E}^.,$$
 o\`u ${\cal D}^\dagger_{{\cal X}\gets {\cal X}^\#}$ 
est obtenu par passage \`a la limite en la variable 
$i$ et $m$ du faisceau: $${\cal D}^{(m)}_{X_i\gets X^\#_i}=
\omega_{X_i}^{-1}\otimes_{{\cal O}_{X_i}}{\cal D}^{(m)}_{X_i}\otimes_{{\cal O}_{X_i}} \omega_{X_i^\#},$$ 
o\`u $(\omega_{X_i})^{-1}$ 
est l'inverse du faisceau inversible $\omega_{X_i}$ des différentielles sur $X_i $ .\\
On a le lemme suivant.  
\begin{surlem}
\begin{enumerate}\label{subsection-lem1}
\item[(i)] Le foncteur $\pi_+$ est un foncteur de 
$D^b_{coh}({\cal D}^\dagger_{{\cal X}^\#,{\bf Q}})$ vers  
 $D^b_{coh}({\cal D}^\dagger_{{\cal X},{\bf Q}})$. 
\item[(ii)] Soit ${\cal E}$ un ${\cal D}^\dagger_{{\cal X}^\#,{\bf Q}}$-module coh\'erent. 
Alors il existe un isomorphisme canonique dans  
$D^b_{coh}(\Ddag_{\XX,\Qr}(^{\dagger}Z))$
$$ \Ddag_{\XX,\Qr}(^{\dagger}Z)\ot_{\Ddag_{\XX^\#,\Qr}}^{\bf L}\EE \simeq
\Ddag_{\XX,\Qr}(^{\dagger}Z)\ot_{\Ddag_{\XX,\Qr}}\pi_+\EE,$$
ces deux complexes étant concentrés en degré $0$.
\end{enumerate}
 \end{surlem} 
\dem. Comme $\Ddag_{\XX^\#,\Qr}$ est de dimension cohomologique finie, il 
suffit, par dévissage, de montrer que $\pi_+(\Ddag_{\XX^\#,\Qr})\in
D^b_{coh}(\Ddag_{\XX,\Qr})$.  
Le faisceau $\Ddag_{\XX^\#\rig \XX,\Qr}$ n'est autre que 
$\Ddag_{\XX,\Qr}$, vu 
comme $\Ddag_{\XX^\#,\Qr}\times \Ddag_{\XX,\Qr}$-bimodule. On a donc 
$$\pi_+\Ddag_{\XX^\#,\Qr}=\omega_{\XX}^{-1}\ot_{\OO_{\XX}}\Ddag_{\XX,\Qr}\ot_{\OO_{\XX}}\omega_{\XX^\#}.$$
Or, comme $\Ddag_{\XX,\Qr}$-module à gauche,
 ce faisceau est isomorphe au 
$\Ddag_{\XX,\Qr}$-module à gauche cohérent et 
induit $$\Ddag_{\XX,\Qr}\ot_{\OO_{\XX}}(\omega_{\XX^\#}\ot_{\OO_{\XX}}\omega_{\XX}^{-1}),$$ 
(où la structure de $\Ddag_{\XX^,\Qr}$-module à gauche vient de 
celle de $\Ddag_{\XX^,\Qr}$). 
En effet, ces différents faisceaux sont isomorphes comme 
$\OO_{\XX}$-modules. Il suffit de vérifier en restriction à 
$\UU$, qu'ils sont isomorphes comme $\Ddag_{\XX,\Qr}$-modules à gauche, 
ce qui provient de l'isomorphisme de 
transposition de 1.3.4 de \cite{Be-smf}. Cela achève le (i).

 Pour la seconde assertion, la fl\`eche $\omega_{\XX}\rig\omega_{\XX^\#} $ induit un 
morphisme apr\`es 
tensorisation par $\OO_{\XX,\Qr}(^{\dagger }Z)$, qui est un 
isomorphisme sur $ \UU$ et 
est donc un isomorphisme puisque l'extension $\OO_{\XX,\Qr}(^{\dagger }Z)\hrig j_*\OO_{{\cal U},\Qr}$ 
est fid\`element plate par \cite{Be1}, 4.3.7.1. 

Il existe un morphisme canonique de $\Ddag_{\XX,\Qr}(^{\dagger}Z)× 
\Ddag_{\XX^\#,\Qr}$-bi-modules 
$$\Ddag_{\XX,\Qr}(^{\dagger}Z)\ot_{\Ddag_{\XX,\Qr}}\DD^{\dagger}_{\XX\lrig 
\XX^{\#}}\rig
\tilde{\omega}_{\XX^{\#}}\ot_{\OO_{\XX}(^{\dagger 
}Z)}\Ddag_{\XX,\Qr}(^{\dagger}Z)\ot_{\OO_{\XX}(^{\dagger }Z)}
\tilde{\omega}_{\XX}^{-1},$$
le deuxi\`eme module \'etant muni de la structure de 
$\Ddag_{\XX^\#,\Qr}$-module
\`a droite provenant de la structure gauche de 
$\Ddag_{\XX,\Qr}(^{\dagger}Z)$ tordue par
$\tilde{\omega}_{\XX^\#}$ et de la structure de 
$\Ddag_{\XX,\Qr}(^{\dagger}Z)$-module \`a
gauche de $\Ddag_{\XX,\Qr}(^{\dagger}Z)$-module, provenant de la 
structure de 
$\Ddag_{\XX,\Qr}(^{\dagger}Z)$-module \`a droite, tordue par 
$\tilde{\omega}_{\XX}^{-1}$. 
Cette fl\`eche est un isomorphisme en restriction \`a $\UU$. Cette 
fl\`eche est donc un isomorphisme de 
$\Ddag_{\XX,\Qr}(^{\dagger}Z)$-modules \`a gauche. On 
v\'erifie en restriction \`a $\UU$, que cet isomorphisme
est aussi $\Ddag_{\XX^\#,\Qr}$-lin\'eaire \`a droite, si bien que 
l'isomorphisme est 
$\Ddag_{\XX^\#,\Qr}$-lin\'eaire \`a droite. En tensorisant \`a droite, 
sur $\Ddag_{\XX^\#,\Qr}$, 
par $\EE$, on en d\'eduit un isomorphisme dans $D^b_{coh}(\Ddag_{\XX,\Qr}(^{\dagger}Z))$ 
$$\Ddag_{\XX,\Qr}(^{\dagger}Z)\ot_{\Ddag_{\XX,\Qr}}
\DD^{\dagger}_{\XX\lrig \XX^{\#}}
\ot_{\Ddag_{\XX^\#,\Qr}}^{\bf L}\EE \simeq 
 \tilde{\omega}_{\XX^{\#}}\ot_{\OO_{\XX}(^{\dagger}Z)}\Ddag_{\XX,\Qr}(^{\dagger}Z)
\ot_{\OO_{\XX}(^{\dagger }Z)}
\tilde{\omega}_{\XX}^{-1}\ot_{\Ddag_{\XX^\#,\Qr}}^{\bf L}\EE.$$
On a une fl\`eche canonique 
$$\Ddag_{\XX,\Qr}(^{\dagger}Z)\rig 
\tilde{\omega}_{\XX^{\#}}\ot_{\OO_{\XX}(^{\dagger }Z)}
\Ddag_{\XX,\Qr}(^{\dagger}Z)\ot_{\OO_{\XX}(^{\dagger }Z)} 
\tilde{\omega}_{\XX}^{-1},$$
envoyant $1$ sur $s\ot 1 \ot s'$, o\`u $s$ est une section locale de 
$\tilde{\omega}_{\XX}$ et $s'$ la section duale de 
$\tilde{\omega}_{\XX^{\#}}^{-1}$. Cette
fl\`eche est un isomorphisme 
bi-$\Ddag_{\XX,\Qr}(^{\dagger}Z)$-lin\'eaire en restriction \`a
$\UU$ (de nouveau grâce à 1.3.4 de \cite{Be-smf}), ce qui donne le fait que c'est un isomorphisme bi-
$\Ddag_{\XX,\Qr}(^{\dagger}Z)\times \Ddag_{\XX^\#,\Qr}$-lin\'eaire et donne la formule du 
(ii) du lemme, puisque $\Ddag_{\XX,\Qr}(^{\dagger}Z)$ est plat sur 
$\Ddag_{\XX,\Qr}$. En restriction à $\UU$, la log-structure devient triviale et 
le faisceau $\pi_+\EE$ s'identifie à $\EE$ comme $\Ddag_{\UU,\Qr}$-module. Pour 
$i\leq -1$, les faisceaux 
de cohomologie $$\Ddag_{\XX,\Qr}(^{\dagger}Z)\ot_{\Ddag_{\XX,\Qr}}\HH^i(\pi_+\EE)$$ sont des $\Ddag_{\XX,\Qr}(^{\dagger}Z) $-modules
cohérents nuls en restriction à $\UU$, donc nuls, ce qui donne le (ii).

\begin{para} Soit $E^\dagger$ un $F$-isocristal 
surconvergent sur $U$. On notera $\tilde{\cal D}^\dagger(.)$ (resp. ${\cal 
D}^\dagger(.)$) le foncteur qui associe \`a $E^\dagger$ sa structure de 
canonique de ${\cal D}^\dagger_{\cal X}(^\dagger {\cal Z})$- 
( resp. ${\cal D}^\dagger_{\cal X})$-module (cf [2],4).  Avec ces notations, nous avons l'isomorphisme de foncteurs: 
 $${\cal D}^\dagger(.)\simeq {\bf res}\circ\tilde{\cal D}^\dagger(.).$$  
Supposons \`a pr\'esent que $E^\dagger$ soit unipotent. Soient $E$ 
le log-isocristal sur $X^\#$ muni d'une structure  de Frobenius 
non-d\'eg\'en\'er\'e lui correspondant d'apr\`es \cite{M-T_unip} et ${\cal 
E}^\dagger={\cal D}^\dagger(E^\dagger)$. Rappelons que la restriction 
de $E$ \`a l'ouvert $U$ (o\`u la log-structure devient triviale)  
n'est autre que le $F$-isocristal convergent associ\'e \`a $E^\dagger$ 
par simple restriction au tube de $U$.  Sous nos hypoth\`eses, $E$ 
correspond de mani\`ere analogue \`a \cite{Be-smf}, 4.6.3 \`a la donn\'ee d'un 
${\cal O}_{{\cal X},{\bf Q}}$-module coh\'erent ${\cal  E}:=sp_*E$ 
muni d'une structure de $F$-${\cal D}^\dagger_{{\cal  X}^\#,{\bf 
Q}}$-module. 

L'énoncé suivant permet de calculer la
cohomologie de $E$ et $E^\dagger$.
\end{para}
 
\begin{surlem}\label{subsection-lem2} 
 Nous conservons les hypoth\`eses et notations pr\'ec\'edentes.  
Soit $f:{\cal X}\to Spf\ W$ le morphisme structural de ${\cal X}/W$. Alors  
\begin{enumerate} 
\item[(i)]  
Pour tout isocristal $E^\dagger$ surconvergent sur $U$, on a 
$$f_+({\cal E}^\dagger)\simeq {\bf R}\Gamma_{rig}(U/K,E^\dagger)[1].$$ 
\item[(ii)]  Pour tout 
cristal $E$ sur $X^\#/W$, on a  $$f_+(\pi_+({\cal E}))\simeq 
{\bf R}\Gamma_{cris}(X^\sharp/W,E)[1]\otimes K$$  
 
\end{enumerate}\end{surlem} 

\dem. La premi\`ere assertion r\'esulte de \cite{Be-survey}, 4.3.6.3  et de
\cite{Be-coh-rig-et-D-mod}, 4.1.7. Passons à la deuxième.
 Remarquons d'abord que $(f_+\circ \pi_+)\EE\simeq f_+^{\#}\EE$. 
En effet, le terme de gauche est égal à 
$$\omega_{\XX}\ot_{\Ddag_{\XX,\Qr}}^{\bf L}{\cal D}^\dagger_{{\cal X}\gets {\cal X}^\#},$$
et donc à 
$$ \omega_{\XX}\ot_{\Ddag_{\XX,\Qr}}{\cal D}^\dagger_{{\cal X}\gets {\cal X}^\#}$$ puisque 
${\cal D}^\dagger_{{\cal X}\gets {\cal X}^\#}$ est un $\Ddag_{\XX,\Qr}$-module à gauche 
plat d'après la démonstration de~\ref{subsection-lem1}, soit encore à 
$\omega_{\XX}\ot_{\Ddag_{\XX,\Qr}}\Ddag_{\XX,\Qr}\ot_{\OO_{\XX}}
(\omega_{\XX^\#}\ot_{\OO_{\XX}}\omega_{\XX}^{-1})$, c'est-à-dire finalement à 
$\omega_{\XX^{\#}}$. On vérifie en restriction à $\UU$ que cet isomorphisme 
est un isomorphisme de $\OO_{\SS}\times \Ddag_{\XX^\#,\Qr}$-bi-modules. Comme le faisceau 
${\cal D}^\dagger_{{\cal S}\gets {\cal X}^\#}$ co\"{\i}ncide avec $\omega_{\XX^{\#}}$
comme $\OO_{\SS}\times \Ddag_{\XX^\#,\Qr}$-bi-module, cela donne l'assertion. Il reste 
à calculer $f_+^{\#}\EE$. Pour cela, on utilise la résolution
$\Ddag_{\XX^\#,\Qr}$-linéaire à droite de $\omega_{\XX^{\#}}$
par le complexe de de Rham suivant dont les termes sont placés en degrés $-1$ et $0$
$$0\rig \Ddag_{\XX^\#,\Qr}\rig \omega_{\XX^{\#}}\ot_{\OO_{\XX}}\Ddag_{\XX^\#,\Qr}\rig 0,$$
et qui est défini localement par $P\mapsto dt/t\otimes \der_{[1]}$, en reprenant les notations de ~\ref{subsection-desc_loc}. 
Remarquons qu'une section locale 
$P$ de $\Ddag_{\XX^\#,\Qr}$ peut se mettre sous la forme 
$P=\sum_{k\in\Ne}\der_{[ k ]}b_k$ tels qu'il existe $C>0, \eta<1$ tels que 
$|b_k|<C\eta^k$. Posons $$Q=\sum_{k\geq 1}\frac{\der_{[ k ]}}{k+1}b_k,$$
dont on vérifie que $Q$ est une section locale de $\Ddag_{\XX^\#,\Qr}$. Alors 
$P$ se décompose $P=b_0+\der_{[1]}Q$, ce qui montre que le complexe précédent 
est quasi-isomorphe à $\omega_{\XX^{\#}}$ placé en degré $0$. Cela permet 
de vérifier que $f_+^{\#}\EE$ s'identifie au complexe de de Rham logarithmique 
de $\EE$ décalé de $1$ et donc à ${\bf R}\Gamma_{cris}(X^\sharp/W,E)[1]\otimes K$,
d'après 4.1.1 de \cite{LS-T_logcris}. 
\begin{surthm} \label{subsection-thm1} 
Soit $E$ un log-isocristal sur $X^\#$ muni d'une  structure de Frobenius 
non-d\'eg\'en\'er\'e et notons $E^\dagger$, le $F$-isocristal surconvergent le long de 
$X\backslash U$ associ\'e \`a $E$ dans \cite{LS-T_logcris}. Alors  ${\cal D}^\dagger(E^\dagger)$ 
a une structure de  $F$-${\cal D}^\dagger_{{\cal X},\Qr}$-module coh\'erent. \end{surthm} 
\dem.  L'isocristal $E$ correspond \`a la donn\'ee d'un ${\cal O}_{{\cal X}_K}$-module \`a 
connexion \`a p\^oles logarithmiques que nous noterons $E_K$ et ${\cal E}:=sp_*E$ \`a celle d'un  
${\cal O}_{{\cal X},{\bf Q}}$-module coh\'erent muni d'une  structure 
de  $F$-${\cal D}^\dagger_{{\cal X}^\#,\Qr}$-module. L'image de 
ce module \`a connexion par le foncteur de \cite{LS-T_logcris} est 
$sp_*(j^\dagger E_K)$. Sa structure de  
${\cal D}^\dagger_{{\cal X},\Qr}(^\dagger Z)$-module coh\'erent est  
$$\tilde{\cal D}^\dagger(E^\dagger)\simeq {\cal D}^\dagger_{\cal X}(^\dagger Z)
\otimes_{{\cal D}^\dagger_{{\cal X}^\#,{\bf Q}}}{\cal E}.$$ 
Pour voir ceci, il suffit d'apr\`es \cite{Be1}, 4.3.12 de 
trouver un homomorphisme de  ${\cal D}^\dagger_{{\cal X},\Qr}(^\dagger Z)$-modules coh\'erents entre ces deux modules qui soit un 
isomorphisme lorsqu'on le restreint \`a $\cal U$. Ce morphisme est 
construit de la mani\`ere suivante: On a une fl\`eche canonique (cf 
\cite{Be-finit}, 1.2.1): $${\cal E}\to \tilde{\cal D}^\dagger(E^\dagger)$$ 
qui induit apr\`es extension des scalaires une fl\`eche 
$${\cal D}^\dagger_{\cal X}(^\dagger Z)\otimes_{{\cal 
D}^\dagger_{{\cal X}^\#}}{\cal E} \to {\cal D}^\dagger_{{\cal 
X}}(^\dagger Z)\otimes_{{\cal D}^\dagger_{{\cal X}^\#,{\bf 
 Q}}}\tilde{\cal D}^\dagger(E^\dagger)$$ o\`u $\tilde{\cal D}^\dagger(E^\dagger)$ 
est consid\'er\'e comme un 
${\cal D}^\dagger_{{\cal X}^\#}$-module. 
Enfin on compose cette fl\`eche avec la fl\`eche canonique 
  $$ {\cal D}^\dagger_{\cal X}(^\dagger Z)\otimes_{{\cal 
D}^\dagger_{{\cal X}^\#}}\tilde{\cal D}^\dagger(E^\dagger)\to \tilde{\cal D}^\dagger(E^\dagger)$$ 
pour obtenir la fl\`eche souhait\'ee. Par le (ii) du lemme
\ref{subsection-lem1}, on a donc $\tilde{\cal D}^\dagger(E^\dagger)\simeq
\Ddag_{\XX,\Qr}(^{\dagger}Z)\ot_{\Ddag_{\XX,\Qr}}\pi_+\EE$ et
on d\'eduit de \cite{Be-survey}, 5.3.6, et en particulier de la
description qui y est donn\'ee du foncteur ${\bf R}j_*j^*$ ($j:U\to X$
est l'immersion ouverte canonique), un triangle distingu\'e 
$${\bf R}\Gamma_{Z}^\dagger(\pi_+{\cal E})\to \pi_+{\cal E}\to {\cal D}^\dagger(E^\dagger),$$ 
o\`u 
${\bf R}\Gamma_{Z}^\dagger$ est le foncteur d\'efini dans \cite{Be-survey}, 
4.4.4. Comme $\pi_+{\cal E}$ est dans $D^b_{coh}({\cal D}^\dagger_{{\cal X},\Qr})$, 
pour montrer que ${\cal D}^\dagger(E^\dagger)$ 
est ${\cal D}^\dagger_{{\cal X},\Qr}$-coh\'erent, il suffit de 
montrer que ${\bf R}\Gamma_{Z}^\dagger(\pi_+{\cal 
E})$ est dans $D^b_{coh}({\cal D}^\dagger_{{\cal X},\Qr})$. On a en fait beaucoup mieux:

\begin{lem} Le complexe ${\bf R}\Gamma_{Z}^\dagger(\pi_+{\cal E})$ est quasi-isomorphe au complexe nul.
\end{lem}

\dem.   En appliquant le foncteur $f_+$ au triangle distingu\'e 
pr\'ec\'edent, on obtient gr\^ace au \ref{subsection-lem2}, le triangle 
distingu\'e $$f_+{\bf R}\Gamma_{Z}^\dagger(\pi_+{\cal E})\to {\bf 
R}\Gamma_{cris}(X^\sharp/W,E)\to {\bf R}\Gamma_{rig}(U/K,E^\dagger).$$   
Comme $E$ est muni d'un Frobenius non-d\'eg\'en\'er\'e, les endomorphismes 
r\'esidus de sa connexion ont pour seule valeur propre 0 (voir \cite{Tr-synt}, 
preuve du th\'eor\`eme 2.5) et donc par \cite{LS-T_logcris}, 4.2 (ii), la deuxi\`eme 
fl\`eche du triangle est un quasi-isomorphisme et ainsi  $f_+{\bf R}\Gamma_{Z}^\dagger(\pi_+{\cal E})$ 
est quasi-isomorphe au complexe 
nul. Comme   ${\bf R}\Gamma_{Z}^\dagger\simeq \oplus_{x\in Z}{\bf R}\Gamma_{x}^\dagger$ 
par \cite{Be-finit}, 2.4 (ii), nous pouvons supposer 
 que $Z$ est un point ferm\'e $u:x\hookrightarrow X$. De plus, 
on a alors par \cite{Be-survey}, 4.4.5 un isomorphisme de foncteurs 
$${\bf R}\Gamma_{Z}^\dagger\simeq u_+u^!.$$ 
Le foncteur compos\'e $(fu)_+$ est le foncteur restriction de la 
cat\'egorie des $K_x$-espaces vectoriels dans celle des  $K$-espaces 
vectoriel, o\`u $K_x=Frac(W(k(x))/K$ est une extension finie. Comme 
$(fu)_+u^!(\pi_+{\cal E})=0$, il en est donc de m\^eme  pour $u^!(\pi_+{\cal E})$ 
et pour $u_+u^!(\pi_+{\cal E})$

On d\'eduit du lemme pr\'ec\'edent que  
$$\pi_+({\cal E})\simeq {\cal D}^\dagger(E^\dagger)$$ et ainsi 
${\cal D}^\dagger(E^\dagger)$ est coh\'erent. De plus, $\pi_+(\EE)$ est concentré en degré 
$0$. 

\begin{surcor} ${\cal D}^\dagger(E^\dagger)$ est holonome.
\end{surcor}

\dem. D'apr\`es [6], 5.3.5, il suffit de montrer qu'il existe un ouvert non-vide de $\cal X$ tel que la restriction 
de ${\cal D}^\dagger(E^\dagger)$ \`a cet ouvert soit le ${\cal D}$-module arithm\'etique associ\'e \`a un $F$-isocristal convergent. 
C'est bien le cas ici puisque, par fonctorialit\'e des constructions la restriction de ${\cal D}^\dagger(E^\dagger)$ \`a $\cal U$ 
n'est autre que le ${\cal D}$-module arithm\'etique associ\'e au $F$-isocristal convergent induit par
 la restriction de $E^\dagger$ au tube de $U$.
\section{Cas g\'en\'eral} 
\bigskip   
 Soient une courbe propre, lisse int\`egre connexe $X/k$, $D$ un diviseur \`a croisement normaux de  $X$ et $U:=X\setminus D$. 
Soit $E^\dagger$ un $F$-isocristal surconvergent sur  $U$. Comme il a \'et\'e d\'emontr\'e dans \cite{M-T_unip}, cet 
isocristal devient unipotent apr\`es un \'eventuel changement de bases par un morphisme fini sur $X$ et \'etale sur $U$. 
Nous allons alors montrer que ${\cal D}^\dagger(E^\dagger)$ est un ${\cal D}^\dagger_{\cal X}$-module 
coh\'erent. Ceci r\'esultera du cas unipotent ainsi que de~\ref{subsubsection-trace}.  

\begin{surprop}\label{subsection-prop2} Soit $E^\dagger$ un $F$-isocristal surconvergent sur $U$. 
Alors ${\cal D}^\dagger(E^\dagger)$ a une  
structure de ${\cal D}^\dagger_{{\cal X},\Qr}$-module holonome.\end{surprop} 
\dem. D'apr\`es \cite{Be-survey}, p.70, exemple 
(i), il suffit de montrer que ${\cal D}^\dagger(E^\dagger)$ 
est ${\cal D}$-coh\'erent. D'apr\`es \cite{M-T_unip}, il existe un 
recouvrement fini et plat $f_0$ : $X'\rig X$, \'etale au-dessus du $U$, tel que 
$E'^\dagger=f_0^{rig,*}(E^\dagger)$ soit unipotent. D'apr\`es \cite{crew-fin}, 
8.3 on peut choisir un rel\`evement ${\cal X}'$ de $X'$ tel que le 
carr\'e commutatif
$$\begin{array}{rcl}  
U'&\rig&U  \\
 \downarrow & &  \downarrow \\ 
X'&\rig&X\end{array}$$
se rel\`eve en un carr\'e
$$\begin{array}{rcl} {\cal U}'&\sta{g}{\rig}&{\cal U} \\
    j'\downarrow& &\downarrow j \\
 {\cal X}'&\sta{f}{\rig}&{\cal X},\end{array}$$ 
o\`u la fl\`eche $f$ est 
finie et plate et la fl\`eche $g$ est bien \'etale puisque c'est le cas de 
la fl\`eche $U'\to U$ par le principe d'invariance topologique des morphismes \'etales (\cite{mi}, 3.23). 
Posons $\EE={\cal D}^\dagger(E^\dagger)$ et $\EE'=f^!\EE$. 
 Comme $f$ est plat, on déduit 
de 2.1.3 
de~\cite{Car_FL} que $\EE'={\bf res}\circ \tilde{f}^!\tilde{{\cal D}}^\dagger(E^\dagger)$. 
On peut aussi déduire de 2.1.6 de~\cite{Car_FL} que 
$f_+(\EE')={\bf res}\circ \tilde{f}_+\tilde{f}^!\tilde{{\cal D}}^\dagger(E^\dagger)$. 
Le module $\tilde{f}^!\tilde{{\cal D}}^\dagger(E^\dagger)$ s'identifie 
à $\tilde{{\cal D}}(f_0^{rig,*}(E^\dagger))$ d'après le résultat rappelé en 
~\ref{subsubsection-comp_iminv}. Le module $\EE'$ est 
donc $\Ddag_{\XX',\Qr}$-cohérent d'après~\ref{subsection-thm1}.

Finalement, le diagramme obtenu en~\ref{subsubsection-trace} donne un diagramme de
$\Ddag_{\XX,\Qr}$-modules
$$\EE\sta{ad}{\rig}{f}_+f^!\EE \sta{tr}{\rig}\EE,$$
et la composée de l'application trace avec l'adjonction est égale à la multiplication 
par $d$, le degré générique de $f$. En particulier, l'application $1/d\cdot tr$ fournit 
un scindage de $ad$ et permet d'identifier $\EE$ à un facteur direct de ${f}_+f^!\EE$, 
comme $\Ddag_{\XX,\Qr}$-module. Le module $f^!\EE$ est $\Ddag_{\XX',\Qr}$-cohérent.  
La cohérence est préservée par $f_+$, car $f$ est propre (4.3.8 de \cite{Be-survey}). 
On en déduit que $\EE$ est un $\Ddag_{\XX,\Qr}$-module cohérent.

\section{Sur la conjecture D de Berthelot}

Nous souhaitons \'etudier le rapport entre la 
proposition pr\'ec\'edente et la conjecture (D) de Berthelot. Rappelons tout d'abord celle-ci:   
\begin{surconj}$\hbox{}$\label{conjectures} 
\item (\cite{Be-survey}, 5.3.6) 
 Si ${\cal E}^\dagger$ est un 
$F$-${\cal D}^\dagger_{{\cal X},\Qr}(^\dagger Z)$-module coh\'erent 
dont la restriction \`a ${\cal U}$ est holonome, alors ${\bf res}({\cal 
E}^\dagger)$ est  
un $F$-${\cal D}^\dagger_{{\cal X},\Qr}$-module 
holonome.
\end{surconj}
Dans \cite{caro}, 4.3.5  la preuve de la conjecture (D) est r\'eduite (dans le cas des courbes) au resultat suivant de Berthelot dont la d\'emonstration devrait apparaitre dans \cite{Be4}:

\begin{prop}\label{caract} Si ${\cal E}^\dagger$ est un $F$-${\cal D}^\dagger_{{\cal X},\Qr}(^\dagger Z)$-module 
coh\'erent dont la restriction \`a ${\cal U}$ est un $F$-isocristal convergent sur $U$, alors ${\cal E}^\dagger$ est 
en fait surconvergent. 
\end{prop} 
 
Nous rappelons au lecteur la d\'emonstration de la conjecture (D) dans le cas des courbes lisse et en profitons pour pr\'eciser l'un des arguments de Caro 
(plus pr\'ecis\'ement son recours \`a \cite{Be1}, 4.3.12). 
\begin{surthm} Nous conservons les hypoth\`eses et notations pr\'ec\'edentes. Soit ${\cal E}^\dagger$ un $F$-${\cal D}^\dagger_{{\cal X},\Qr}(^\dagger Z)$-module 
coh\'erent dont la restriction \`a ${\cal U}$ est holonome. 
Alors ${\cal E}^\dagger$ un $F$-${\cal D}^\dagger_{{\cal X},\Qr}$-module holonome.
\end{surthm}

\dem. Comme la restriction \`a ${\cal U}$ de ${\cal E}^\dagger$ est holonome, 
il existe d'apr\`es \cite{Be-survey}, p.70 
un ouvert non-vide ${\cal V}$ de $\cal U$ tel que ${\cal E}^\dagger|_{{\cal V}}$ corresponde \`a un $F$-isocristal
convergent. Nous pouvons \'ecrire la fibre sp\'eciale $V$ de ${\cal  V}$ sous la forme $V=X\setminus (Z\cup Z')$, 
pour un certain diviseur $Z'$ disjoint de $Z$. En particulier, ${\cal E}^\dagger|_{{\cal V}}$ 
est ${\cal O}_{{\cal V},{\bf Q}}$-coh\'erent et donc gr\^ace \`a la caract\'erisation des isocristaux 
surconvergents de \ref{caract}, on en d\'eduit que le 
$F$-${\cal D}^\dagger_{{\cal X},\Qr}(^\dagger Z\cup Z')$-module coh\'erent 
${\cal E}^\dagger={\cal  E}^\dagger\otimes_{{\cal D}^\dagger_{\cal X,\Qr}(^\dagger Z)}
 {\cal D}^\dagger_{\cal X,\Qr}(^\dagger Z\cup Z')$ correspond \`a un F-isocristal surconvergent
 sur $V/K$. Par la proposition \ref{subsection-prop2}, ${\bf res}({\cal E}'^\dagger)$ a ainsi une structure 
de ${\cal  D}^\dagger_{\cal X,\Qr}$-module holonome. On peut factoriser le foncteur ${\bf res}(.)$ en 
le foncteur compos\'e:
$${\bf res}(res_{Z\cup Z',Z'})(.),$$
o\`u le foncteur $res_{Z\cup Z',Z'}$ consiste \`a consid\'erer tout 
${\cal D}^\dagger_{\cal X,\Qr}(^\dagger Z\cup Z')$ comme un ${\cal D}^\dagger_{\cal X,\Qr}(^\dagger Z')$-module. 
On peut montrer que  
$$res_{Z\cup Z',Z'}({\cal E}'^\dagger)=
{\bf res}({\cal E}^\dagger)\otimes_{{\cal D}^\dagger_{\cal X,\Qr}} {\cal D}^\dagger_{\cal X,\Qr}(^\dagger Z').$$ 
Pour montrer cette assertion il suffit en effet de regarder ce qui se 
passe sur le complementaire de $Z'$, not\'e $U'$. L'assertion devient alors
$${\bf res}({\cal E}^\dagger)=
{\bf res}({\cal E}^\dagger\otimes_{{\cal D}^\dagger_{\cal X,\Qr}
(^\dagger Z)} {\cal D}^\dagger_{\cal X,\Qr}(^\dagger Z))$$ 
au-dessus de l'ouvert $U'$. Et donc, on conclut en remarquant que 
${\cal E}^\dagger={\cal E}^\dagger\otimes_{{\cal D}^\dagger_{\cal X,\Qr}(^\dagger Z)} 
{\cal D}^\dagger_{\cal X,\Qr}(^\dagger Z)$ d'apr\`es \cite{caro3}, 2.2.8. On a d'autre part un triangle distingu\'e
$${\bf R}\Gamma^\dagger_{Z'}({\bf res}({\cal E}^\dagger))\to {\bf res}({\cal 
E}^\dagger)\to {\bf res}({\bf res}({\cal E}^\dagger)\otimes_{{\cal D}^\dagger_{\cal X,\Qr}}
 {\cal D}^\dagger_{\cal X,\Qr}(^\dagger Z')),$$
o\`u le terme de droite est holonome puisqu'il n'est autre que
${\bf res}({\cal E}'^\dagger)$ par l'assertion pr\'ec\'edente. Il
reste donc \`a v\'erifier l'holonomie de ${\bf  R}\Gamma^\dagger_{Z'}({\bf res}({\cal E}^\dagger))$. 
Notons $u':Z'\hookrightarrow X$, et $u'':Z'\hookrightarrow U$ les immersions ferm\'ees canoniques.
 D'apr\`es \cite{Be-survey}, 4.4.5, on a 
$${\bf  R}\Gamma^\dagger_{Z'}({\bf res}({\cal E}^\dagger))=u'_+u'^!({\bf res}({\cal
  E}^\dagger)).$$
Comme $Z'$ est disjoint de $Z$, on a $u'^!({\bf res}({\cal  E}^\dagger))=u''^!{\bf res}({\cal  E}^\dagger)|_{\cal U}$ 
qui a une structure de ${\cal  D}^\dagger_{{\cal Z}'}$-module holonome d'apr\`es \cite{Be-survey}, p.72, exemples (iv). 
Finalement, on conclut que ${\bf  R}\Gamma^\dagger_{Z'}({\bf res}({\cal E}^\dagger))$ 
est un ${\cal  D}^\dagger_{\cal X,\Qr}$-module holonome en utilisant 
\cite{Be-survey}, p.72, exemples (iii).  
\bigskip
\begin{surcor}
Soit ${\cal X}/\SS$ un sch\'ema formel lisse de dimension relative
un. Alors la cat\'egorie des 
$F$-${\cal  D}^\dagger_{{\cal X},\Qr}$-complexes holonomes born\'es est stable par les six op\'erations 
cohomologiques de Grothendieck.
\end{surcor}

\dem. D'apr\`es \cite{Be-survey}, il suffit de d\'emontrer les
conjectures A et D (respectivement \cite{Be-survey}, 5.3.6 A) et D)). 
L'assertion r\'esulte alors du th\'eor\`eme
pr\'ec\'edent et de \cite{caro2} o\`u il est d\'emontr\'e que la
conjecture D implique la conjecture A.

\begin{para} Pour terminer, nous souhaitons signaler que les r\'esultats pr\'ec\'edents admettent des \'enonc\'es analogues dans le cas simple o\`u 
les $F$-isocristaux surconvergents correspondent \`a des modules sur des 
anneaux de s\'eries formelles munis d'une connexion et d'un op\'erateur 
Frobenius horizontal. Nous ne faisons ici qu'\'ebaucher ces analogues locaux. Pour plus de pr\'ecision, 
nous renvoyons le lecteur \`a \cite{crew-excep} ou encore \`a \cite{Ma}.
\end{para}

\begin{para} Soit $k$ un corps parfait de caract\'eristique $p$. On 
note $W$ l'anneau des vecteurs de Witt de $k$ et $K$ le corps des fractions de 
$W$. On note $\sigma$ le Frobenius canonique de $K$. Pour une ind\'etermin\'ee 
$x$, on note :

\bigskip
\noindent $B=\lbrace a=\sum_{-\infty}^{+\infty} a_ix^i|\
a_i\in K,\  \sup_i|a_i|<\infty,\ |a_i|\rightarrow 0\ (i\rightarrow 
-\infty)\rbrace.$\\

\noindent $B^\dagger=\lbrace a\in B |\ 
|a_i|r^i\rightarrow 0\ (i\rightarrow -\infty)\hbox{ pour un certain } 
0<r<1 \rbrace.$
\bigskip

Notons indiff\'eremment $R$ l'anneau $B$ ou $B^\dagger$. C'est un corps 
de valuation 
discr\^ete de corps r\'esiduel $k((x))$ muni d'un op\'erateur Frobenius 
et d'une $K$-d\'erivation
$$d: R\rightarrow \omega_{R}:=R.{dx\over x}.$$
Rappelons (\cite{M-T_unip}) qu'un $(\phi,\nabla)$-module sur $R$ 
consiste en
la donn\'ee d'un $R$-module libre de rang fini $M$ 
muni d'une connexion $K$-lin\'eaire :
$$\nabla : M\rightarrow M\otimes \omega_{R}$$
telle que $\nabla(am)=da\otimes m+a\nabla(m)$ pour tout $a\in R$ 
et tout $m\in M$. Le module est d'autre part muni d'un 
endomorphisme $\sigma$-lin\'eaire :
$$\phi: M\rightarrow M$$
tel que $\phi(M)$ engendre  $M$ et :
$$\nabla \phi=(\sigma\otimes \phi)\nabla.$$

En fait, les $(\phi,\nabla)$-modules sur $B$ (resp. $B^\dagger$) 
correspondent moralement aux $F$-isocristaux convergents (resp. 
surconvergents) sur $k((x))$.\\
\end{para}

\begin{para} Posons ${\cal X}:=Spf (W[[x]])$ et ${\cal U}=Spf(\widehat{W[[x]][{1\over 
x}]})$ et $D^\dagger_{\cal X}$ et $D^\dagger_{\cal U}$ (resp. 
$D^\dagger(0)$) l'anneau  des op\'erateurs diff\'erentiels surconvergents sur 
${\cal X}$ et ${\cal U}$ (resp. l'anneau des op\'erateurs diff\'erentiels 
surconvergents le long du diviseur $(x=0))$. On note $\partial^{[i]}:={1\over i!}{d^i\over dx^i}$ et on 
munit l'anneau des series formelles de la norme
de Gauss $$|\sum_{i\geq 0} a_ix^i|=p^{-min_i(v(a_i))}.$$ Les anneaux differentiels $D^\dagger_{{\cal X},\Qr}$, 
$D^\dagger_{{\cal U},\Qr}$ et $D^\dagger(0)_\Qr$ sont d\'efinis de la mani\`ere suivante:

$$D^\dagger_{{\cal X},\Qr}:=\left\{\sum_{i\geq
0} b_i\partial^{[i]}/ b_i\in W[[x]][{1\over p}]/\exists
c,\eta\in \Rr, |b_i|<c\eta^i\right\},$$
$$D^\dagger_{{\cal U},\Qr}:=\left\{\sum_{i\geq
0} b_i\partial^{[i]}/ b_i\in\widehat{W[[x]][{1\over
x}]}\left[{1\over p}\right]=B/\exists c,\eta\in \Rr,
|b_i|<c\eta^i\right\}$$
$$D^\dagger(0)_\Qr:=\left\{\sum_{i\geq 0}
b_i\partial^{[i]}/ b_i\in B^\dagger/\exists c,\eta\in \Rr, \eta<1, |b_i|<c\eta^i\right\}.$$

Il semble raisonnable de penser qu'on a une \'equivalence de 
cat\'egories entre les $(\phi,\nabla)$-modules sur $B$ et les $B$-espaces 
vectoriels de dimension finie munis d'une structure de 
$F$-$D^\dagger_{\cal U}$-modules (on a $B={\cal O}_{{\cal U},{\bf Q}}$).\\
\end{para} 

La conjecture (D) et la caracterisation des F-isocristaux surconvergents de Berthelot peuvent \^etre traduit de la maniere suivante (et d\'emontrer de mani\`ere analogue au cas des courbes):
\begin{prop}$\hbox{}$\label{conjlocal}
\begin{enumerate}
\item Soit $M(0)$ un $D^\dagger(0)_\Qr$-module coh\'erent tel que $M:=M(0)\otimes_{D^\dagger(0)}
  D^\dagger_{{\cal U}}$ soit un $(\phi,\nabla)$-module sur ${\cal O}_{{\cal U},{\bf Q}}$. Alors ${\bf res}(M)$ est un
  $D^\dagger_{{\cal X},\Qr}$-module coh\'erent.
\item  Soit $M(0)$ un $D^\dagger(0)_\Qr$-module coh\'erent tel que 
$M:=M(0)\otimes_{D^\dagger(0)} D^\dagger_{\cal U}$ soit un 
$(\phi,\nabla)$-module sur $B$. Alors il existe un $(\phi,\nabla)$-module $M^\dagger$ sur 
$B^\dagger$ tel que 
$$M^\dagger\otimes_{B^\dagger} B\simeq M\otimes_{D^\dagger(0)} 
D^\dagger_{\cal U}.$$
\end{enumerate} 
\end{prop}

On peut \'egalement montrer de mani\`ere analogue \`a 
\ref{subsection-prop2}, le r\'esultat local suivant (voir aussi \cite{crew-excep}, 
3.3):
\begin{surprop}\label{prop2local}
Soit $M$ un $(\phi,\nabla)$-modules sur $B^\dagger$. Alors $M$ est 
holonome, i.e. ${\bf res}(M\otimes_{D^\dagger_{\cal X}}D^\dagger(0))$ 
est un $(\phi,\nabla)$-module sur $B^\dagger$.
\end{surprop}

%
\begin{flushleft}
 Christine Noot-Huyghe \\
 Institut de Recherche Math\'ematique Avanc\'ee (IRMA) \\
 Universit\'e Louis Pasteur et CNRS\\
 7 rue Ren\'e Descartes \\
  67084 Strasbourg Cedex France\\
 m\'el huyghe@math-u.strasbg.fr, http://www-irma.u-strasbg.fr/\textasciitilde huyghe
\end{flushleft}
%
%
\begin{flushleft}
  Fabien Trihan\\
  Universit\'e de Mons-Hainaut\\
  « Le Pentagone », 6 avenue du champ de Mars\\
  B-7000 Mons (Belgique)\\
 m\'el trihan@umh.ac.be
\end{flushleft}


\begin{thebibliography}{27}

\bibitem{andre_monp}
Yves Andr{\'e}.
\newblock Filtrations de type {H}asse-{A}rf et monodromie {$p$}-adique.
\newblock {\em Invent. Math.}, 148(2), p.~285--317, 2002.


\bibitem{Be-coh-rig-et-D-mod}
P.~Berthelot.
\newblock { Cohomologie rigide et th\'eorie des $\DD$-modules}.
\newblock {\em { Proc. Conf. $p$-adic Analysis (Trento 1989), Lecture Notes in
  Math., Springer-Verlag}}, 1454, p.~78--124, 1990.

\bibitem{Be1}
P.~Berthelot.
\newblock { $\DD$-modules arithm\'etiques I. {O}p\'erateurs diff\'erentiels de
  niveau fini}.
\newblock {\em Ann.~{scient}.~{\'E}c.~Norm.~Sup.}, $4^ e$ s\'erie, t.~29,
  p.185--272, 1996.

\bibitem{Be-finit}
P.~Berthelot.
\newblock Finitude et puret\'e cohomologique en cohomologie rigide.
\newblock {\em Invent.~Math}, 128, p.~329--377, 1997.


\bibitem{Be-smf}
P.~Berthelot.
\newblock {D-modules arithm\'etiques II descente par Frobenius}.
\newblock {\em Bull. Soc. Math. France}, M\'emoire 81, p.~1--135, 2000.

\bibitem{Be-survey}
P.~Berthelot.
\newblock { Introduction \`a la th\'eorie arithm\'etique des $\cal D$-modules}.
\newblock {\em {Bull. Soc. Math. France}}, Cohomologies $p$-adiques et
  applications arithm\'etiques II, Ast\'erisque No. 279, 2002.

\bibitem{Be4}
P.~Berthelot.
\newblock{$\DD$-modules arithmetiques IV, Vari\'et\'es caract\'eristiques}
\newblock{En pr\'eparation.}

\bibitem{Car_FL}
D.~Caro.
\newblock {Fonctions $L$ associ\'ees aux $\DD$-modules arithm\'etiques}.
\newblock {\em Th\`ese de Doctorat, Universit\'e de Rennes I}, 2002.

\bibitem{caro}
D.~Caro.
\newblock {Fonctions $L$ associ\'ees aux $\DD$-modules arithm\'etiques. Cas des courbes}.
\newblock {Compos. Math. 142, No. 1, 169-206 (2006)}.

\bibitem{caro2}
D.~Caro.
\newblock {D-modules arithmétiques surcohérents. Application aux
  fonctions L}.
\newblock{\em  Ann. Inst. Fourier, Grenoble 54, 6 (2004), 1943-1996.}.

\bibitem{caro3}
D.~Caro.
\newblock {Comparaison des foncteurs duaux des isocristaux 
surconvergents}.
\newblock{\em Rend. Sem. Mat. Univ. Padova 114 
(2005), 81 p.}. 

\bibitem{Cr-mon_Fis}
R.~Crew.
\newblock $F$-isocrystals and their monodromy groups.
\newblock {\em {Ann. Sci. Ecole Norm. Sup. }}, (4) 25 no. 4 , p.~429--464,
  1992.

\bibitem{crew-fin}
R.~Crew.
\newblock {Finiteness theorems for the cohomology of an overconvergent isocrystal on a curve}.  
\newblock {\em {Ann. Sci. Ecole Norm. Sup. }},(4) 31 no. 6 , p.~717--463,
  1998.
\bibitem{crew-excep}
R.~Crew.
\newblock {Arithmetic $\cal D$-modules on a formal curve}. 
\newblock {Math. Ann. 336, No. 2, 439-448 (2006)}.

\bibitem{huy}
C.~Huyghe.
\newblock Construction et \'etude de la transform\'ee de Fourier pour les $\cal D$-modules arithm\'etiques.
\newblock {\em { th\`ese de Doctorat, Universit\'e de Rennes1, 1995.}}

\bibitem{HuFour}
C.~Huyghe.
\newblock {\it $D^{\dagger}$-affinit\'e des sch\'emas projectifs}.
\newblock {\em Ann. Inst. Fourier}, {\bf t.~48}, fascicule 4, p.~913--956,
  1995.

\bibitem{Kato_log1}
K.~Kato.
\newblock Logarithmic structures of Fontaine-Illusie.
\newblock {\em { Algebraic analysis, geometry, and number theory (Baltimore,
  MD, 1989), Johns Hopkins Univ. Press}}, p.~191--224, 1988.

\bibitem{Kedlaya_semi_stable}
Kiran~S. Kedlaya.
\newblock Semistable reduction for overconvergent {$F$}-isocrystals on a curve.
\newblock {\em Math. Res. Lett.}, 10(2-3), p.~151--159, 2003.

\bibitem{kedlaya_monp}
Kiran~S. Kedlaya.
\newblock A {$p$}-adic local monodromy theorem.
\newblock {\em Ann. of Math. (2)}, 160(1), p.~93--184, 2004.

\bibitem{LS-T_logcris}
B.~Le Stum and F.~Trihan.
\newblock Log-cristaux et surconvergence.
\newblock {\em { Ann. Inst. Fourier}}, 51, p.~1189--1207, 2001.

\bibitem{Ma}
A.~Marmora
\newblock Constantes locales $p$-adiques.
\newblock Th\`ese d'Etat, Universit\'e Paris 13, 2006.

\bibitem{M-T_unip}
B.~Matsuda and F.~Trihan.
\newblock Image directe sup\'erieure et unipotence.
\newblock {\em {J. Reine Angew. Math. 569, 47-54 (2004)}}.

\bibitem{mebkhout_monp}
Z.~Mebkhout.
\newblock Analogue {$p$}-adique du th\'eor\`eme de {T}urrittin et le
  th\'eor\`eme de la monodromie {$p$}-adique.
\newblock {\em Invent. Math.}, 148(2), p.~319--351, 2002.

\bibitem{mi}
J.S..~Milne.
\newblock {Etale cohomology.}
\newblock {\em Princeton Mathematical Series}, 33. Princeton University Press, Princeton, N.J., 1980. xiii+323 pp.

\bibitem{Monta_th}
C.~Montagnon.
\newblock {G\'en\'eralisation de la th\'eorie arithm\'etique des $\DD$-modules
  \`a la g\'eom\'etrie logarithmique}.
\newblock {\em Th\`ese de Doctorat, Universit\'e de Rennes I}, 2002.

\bibitem{Tr-synt}
F.~Trihan.
\newblock {Cohomologie syntomique des $F-T$-cristaux.}
\newblock {\em Rend. Sem. Mat. Univ. Padova}, 108, p.~1--26, 2002.

\bibitem{Ts1}
N.~Tsuzuki.
\newblock {Morphisms of $F$-isocrystals and the finite monodromy theorem for
  unit-root $F$-isocrystals}.
\newblock {\em Duke Math. J.}, 111, no. 3, p.~385--418, 2002.

\bibitem{vir}
A.~Virrion.
\newblock{Trace et dualit\'e relative pour les $\cal D$-modules arithm\'etiques}.
\newblock {\em Adolphson, Alan (ed.) et al., Geometric aspects of Dwork theory. Vol. I, II. Berlin: Walter de Gruyter. 1039-1112 (2004)}.

\end{thebibliography}

\end{document}